\documentclass{article}[12pt]

\usepackage{graphicx,tikz,color}
\usepackage{amsmath,amsfonts,amssymb,latexsym}

\newtheorem{theorem}{Theorem}[section]

\newtheorem{proposition}[theorem]{Proposition}

\newtheorem{lemma}[theorem]{Lemma}

\newtheorem{problem}[theorem]{Problem}

\newcommand{\proof}{\noindent{\bf Proof.\ }}
\newcommand{\qed}{\hfill $\square$ \bigskip}
\newcommand{\qqed}{($\Box$)}

\newcommand{\es}{{\rm es}_{\chi}}
\newcommand{\KD}{K_2^{(4)}}
\newcommand{\KT}{K_3^{(2,2,1)}}
\newcommand{\CF}{C_4^{(2,1,2,1)}}
\newcommand{\KKD}{{{\cal F}_{(k,2)}}}
\newcommand{\KCH}{{{\cal F}_{{\chi},{k}}}}

\begin{document}

\title{Critical graphs for the chromatic edge-stability number}

\author{Bo\v stjan Bre\v sar$^{a,b}$\thanks{Email: \texttt{bostjan.bresar@um.si}}
\and Sandi Klav\v zar$^{c,a,b}$\thanks{Email: \texttt{sandi.klavzar@fmf.uni-lj.si}} \and Nazanin Movarraei$^d$\thanks{Email: \texttt{nazanin.movarraei@gmail.com}}
}

\maketitle

\begin{center}
$^a$ Faculty of Natural Sciences and Mathematics, University of Maribor, Slovenia\\
\medskip

$^b$ Institute of Mathematics, Physics and Mechanics, Ljubljana, Slovenia\\
\medskip

$^c$ Faculty of Mathematics and Physics, University of Ljubljana, Slovenia\\
\medskip

$^d$ Department of Mathematics, Yazd University, Iran

\end{center}

\begin{abstract}
The chromatic edge-stability number ${\rm es}_{\chi}(G)$ of a graph $G$ is the minimum number of edges whose removal results in a spanning subgraph $G'$ with $\chi(G')=\chi(G)-1$. Edge-stability critical graphs are introduced as the graphs $G$ with the property that ${\rm es}_{\chi}(G-e) < {\rm es}_{\chi}(G)$ holds for every edge $e\in E(G)$. If $G$ is an edge-stability critical graph with $\chi(G)=k$ and ${\rm es}_{\chi}(G)=\ell$, then $G$ is $(k,\ell)$-critical. Graphs which are $(3,2)$-critical and contain at most four odd cycles are classified. It is also proved that the problem of deciding whether a graph $G$ has $\chi(G)=k$ and is critical for the chromatic number can be reduced in polynomial time to the problem of deciding whether a graph is $(k,2)$-critical. 
\end{abstract}

\noindent
{\bf Keywords:} chromatic edge-stability; edge-stability critical graph; odd cycle; computational complexity \\

\noindent
{\bf AMS Subj.\ Class.\ (2010)}: 05C15, 05C38, 68Q25

\section{Introduction}

Given a graph $G$, a function $c:V(G)\rightarrow [k]=\{1,\ldots, k\}$ such that $c(u)\neq c(v)$ if $uv\in E(G)$, is called a {\em $k$-coloring} of $G$. The minimum $k$ for which $G$ is $k$-colorable is the {\em chromatic number} of $G$, and is denoted by $\chi(G)$. The {\em chromatic edge-stability number}, $\es(G)$, of $G$ is the minimum number of edges of $G$ whose removal results in a graph $G_{1}$ with $\chi(G_{1})=\chi(G)-1$. 

The chromatic edge-stability number was introduced in~\cite{staton-1980}, where $\es$ was bounded from the above for regular graphs in terms of the size of a given graph. Somehow surprisingly, this natural coloring concept only recently received some of the deserved attention. In~\cite{kemnitz-2018} the edge-stability number was compared with the chromatic bondage number and bounded for several graph operations. In~\cite{akbari-2019+}, among other results, several bounds on $\es$ are proved and a Nordhaus-Gaddum type inequality derived. We also mention here the concept of edge-transversal number alias 
bipartite edge frustration, defined as the smallest number of edges that have to be deleted from a graph to obtain a bipartite spanning subgraph, see~\cite{ber-2000, doslic-2007, dvorak-2012, kral-2004, yarahmadi-2011, zhu-2009}. In particular, if $\chi(G) = 3$, then the two invariants coincide. A related study is concerned with the minimum number of edges that an $n$-vertex graph must have so that one can reduce it to a bipartite graph by the removal of a fixed number of edges~\cite{kostochka-2015}. 
 
We say that a graph $G$ is {\em edge-stability critical} if $es_{\chi}(G-e)<es_{\chi}(G)$ holds for every edge $e\in E(G)$. To simplify the writing in this paper, we say that a graph $G$ is {\em $(k,\ell)$-critical}, where $k,\ell\ge 2$, if $G$ is an edge-stability critical graph with $\chi(G)=k$ and $\es(G)=\ell$. 
A graph $G$ is $(k,2)$-critical if and only if for every edge $e$ we have $\chi(G-e)=\chi(G)=k$ and there exists an edge $e'\in E(G)$ such that $\chi(G-\{e,e'\})=k-1$. From this point of view we recall that a {\em double-critical graph} is a connected graph in which the deletion of any pair of adjacent vertices decreases the chromatic number by two. The Erd\H{o}s-Lov\'{a}sz Tihany conjecture asserts that $K_k$ is the only double-critical $k$-chromatic graph~\cite{erdos-1968}. For recent results on this problem see~\cite{roso-2017, stiebitz-2017}. 

Note that $(2,2)$-critical graphs are precisely the graphs with two edges. Hence, if we restrict to isolate-free graphs, then there are only two $(2,2)$-critical graphs, $P_3$ and $2K_2$. Since isolated vertices play no role in this study, we assume that all graphs considered in the rest of the paper are isolate-free. 

To formulate our main result, we introduce the following four families of graphs, where $G + H$ denotes the disjoint union of graphs $G$ and $H$.  Let ${\cal A}=\{C_{2k+1}+ C_{2\ell+1}\,|\,k,\ell\ge 1\}$ and let $\cal B$ be the family of graphs that are obtained from $C_{2k+1} + C_{2\ell+1}$, $k, \ell\ge 1$, by identifying a vertex of $C_{2k+1}$ with a vertex of $C_{2\ell+1}$. Let $x_i,y_i$ be the endvertices of the paths $Q_i$, $i\in[4]$,
exactly two of the $Q_i$ are odd, and at most one of them is of length $1$. The family $\cal C$ consists of the graphs that are obtained from such four paths, by identifying the vertices $x_1$, $x_2$, $x_3$, and $x_4$ and also identifying the vertices $y_1$, $y_2$, $y_3$, and $y_4$. The family $\cal D$ consist of the following subdivisions of the graph $K_4$: (i) all the subdivided paths are of odd length, (ii) exactly three of the paths are odd, and these three paths induce an odd cycle, and (iii) exactly two of the paths are odd, and these two paths are vertex disjoint. Our main result now reads as follows. 

\begin{theorem}
\label{thm:main}
$\cal A\cup \cal B\cup \cal C\cup \cal D$ is the family of $(3,2)$-critical graphs (without isolated vertices) that contain at most four odd cycles.
\end{theorem}

In the next section we prove Theorem~\ref{thm:main}. In Section~\ref{sec:complexity}, we prove that the problem of deciding whether a graph $G$ is critical for the chromatic number and $\chi(G)=k$ can be reduced in polynomial time to the problem of deciding whether a graph is $(k,2)$-critical. We end the paper with concluding remarks concerning a possible  characterization of $(3,2)$-critical graphs. In particular, we find one more family of $(3,2)$-critical graphs, and pose a problem whether the discovered families contain all $(3,2)$-critical graphs.

\section{Proof of Theorem~\ref{thm:main}}

It is clear that graphs with only one odd cycle cannot be $(3,2)$-critical. To describe $(3,2)$-critical graphs with exactly two odd cycles, we first show the following general fact about such graphs. 

\begin{lemma}
\label{lem:general-2-cycles}
If $G$ is a graph that has exactly two odd cycles, and these cycles share no edge, then the cycles intersect either in a single vertex or not at all. 
\end{lemma}

\proof
Let $C$ and $D$ be the odd cycles of $G$ and assume that they have more than one common vertex. Let $u$ and $v$ be vertices from $V(C)\cap V(D)$ such that $d_C(u,v)$ is smallest possible. 
Let $P$ be a $u,v$-subpath on $C$ of length $d_C(u,v)$, and let $P'$ be the other $u,v$-path on $C$. Let $Q$ be a $u,v$-path on $D$, which is internally disjoint from $C$. If the length of $Q$ is of the same parity as the length of $P$, then $Q\cup P'$ is an odd cycle, different from $C$ and $D$. Otherwise, if the length of $Q$ is of different parity as the length of $P$, then $Q\cup P$ is a third odd cycle in $G$, again a contradiction. 
\qed

\begin{theorem}
\label{thm:two-odd-cycles}
The $(3,2)$-critical graphs that contain exactly two odd cycles are precisely the graphs of the families $\cal A$ and $\cal B$. 
\end{theorem}

\proof
Suppose that $G$ contains exactly two odd cycles. If the  two odd cycles have a common edge, then $\es(G)=1$, hence $G$ is not $(3,2)$-critical. Thus, by Lemma~\ref{lem:general-2-cycles}, the two cycles are either disjoint or intersect in a single vertex. There are no other edges in $G$ but the edges of these two cycles, for otherwise, removing an edge $e$ that is not in any of these two cycles gives $\es(G-e)\ge 2$, a contradiction. 
\qed

\begin{lemma}
\label{lem:general}
If $G$ is a graph that has exactly three odd cycles, then the intersection of every two odd cycles is either empty or is a path.   
\end{lemma}
\proof
Let $C=v_1,\ldots,v_k,v_1$ ($k$ odd) be one of the three odd cycles in $G$. 
Suppose that $D$ is an odd cycle in $G$ such that the intersection $C\cap D$ is neither empty nor a path. This implies that  $C\cap D$ is a union of at least two paths on $C$. We may assume without loss of generality that one of these paths is induced by vertices $v_1,\ldots, v_t$, where $t\in [k-1]$. Hence, there is a path $P$ between $v_t$ and a vertex $v_r$ in $V(C)\setminus \{v_1,\ldots,v_t\}$, which is internally in $D\setminus C$ (that is, except for its endvertices $v_t$ and $v_r$, all edges and eventual other vertices of $P$ are not in $C$). Analogously, there is a path $P'$ between $v_1$ and a vertex $v_s$ of $V(C)\setminus \{v_1,\ldots,v_t\}$, which is internally in $D\setminus C$ (that is, except for its endvertices $v_1$ and $v_s$ all edges and eventual other vertices of $P'$ are not in $C$). Note that $P$ and $P'$ are either disjoint or their intersection consists only of vertex $v_1$ (if $t=1$) or vertex $v_r$ (if $r=s$) or both. 

Let $Q$ and $Q'$ be the two paths on $C$ between $v_t$ and $v_r$. If $P$ has the same parity as $Q$, then $P\cup Q'$ is an odd cycle. Otherwise, $P\cup Q$ is an odd cycle. 
Note that any of these cycles is distinct from $C$ and $D$. 
Now, an analogous argument for $P'$ implies that there is an odd cycle $F$ that involves the path $P'$ and a subpath of $C$ between $v_1$ and $v_s$. Clearly, since $P$ and $P'$ are distinct, also $F$ is different from $P\cup Q$ and $P\cup Q'$. We derive that there are at least four odd cycles in $G$, a contradiction. 
\qed

\begin{theorem}
\label{thm:threeoddcylces}
There are no $(3,2)$-critical graphs that contain exactly three odd cycles.
\end{theorem}
\proof
Suppose that there is $(3,2)$-critical graph $G$  that contains exactly three odd cycles, $C$, $D$, and $E$. By Lemma~\ref{lem:general}, the intersections of pairs of cycles are either empty or they are paths. If one of the intersections, say $C\cap D$, has no edges, then removing an edge of $E$, which does not belong to $C\cup D$, yields a graph $G'$ with $\es(G')>1$, hence $G$ is not $(3,2)$-critical. Also, if there exists an edge $e$ that belongs to all three odd cycles, we infer that $\chi(G-e)=2$, again a contradiction. From these two observations we infer that $C\cap D$ is a path $P$, that $C\cap E$ is a path $P'$, and that $D\cap E$ is a path $P''$, where the three paths are non-trivial and are pairwise edge disjoint. The cycle $E$ is not equal to the cycle $C\cdot D$, which is the cycle obtained from $C\cup D$ by removing the internal vertices of $P$, since $C\cdot D$ is an even cycle. 
From this we derive that there exists a path $Q$ from a vertex $x$ in $C\setminus D$ to a vertex $y$ in $D\setminus C$, all of which edges and all eventual internal vertices are not in $C\cup D$. Note that the union of the cycles $C$ and $D$ and the path $Q$ yields a subgraph of $G$, which is a subdivision of $K_4$.

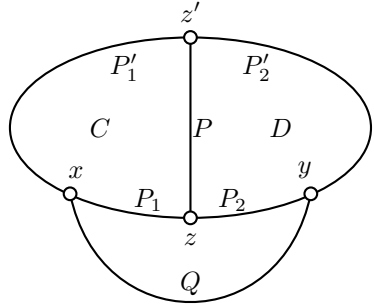
\begin{figure}[ht!]
\begin{center}
\begin{tikzpicture}[scale=0.8,style=thick]
\def\vr{3pt}
\path (0,1) coordinate (z);
\path (0,4) coordinate (z');
\path (-2.0,1.4) coordinate (x);
\path (2,1.4) coordinate (y);
\draw (z) -- (z');
\draw (z) arc(270:90:3cm and 1.5cm);
\draw (z) arc(270:90:-3cm and 1.5cm);
\draw (x) .. controls (-1.5,-1) and (1.5,-1) .. (y);

\draw (z)  [fill=white] circle (\vr);
\draw (z')  [fill=white] circle (\vr);
\draw (x)  [fill=white] circle (\vr);
\draw (y)  [fill=white] circle (\vr);
\draw [below] (z)++(0,-0.1) node {$z$};
\draw [above] (z')++(0,0.1) node {$z'$};
\draw [above] (x)++(0.1,0.1) node {$x$};
\draw [above] (y)++(-0.1,0.1) node {$y$};
\draw (0.2,2.5) node {$P$};
\draw (-1.5,2.5) node {$C$};
\draw (1.5,2.5) node {$D$};
\draw (0,-0.1) node {$Q$};
\draw (-0.7,1.3) node {$P_1$};
\draw (0.7,1.3) node {$P_2$};
\draw (-1.1,3.5) node {$P_1'$};
\draw (1.1,3.5) node {$P_2'$};
\end{tikzpicture}
\end{center}
\vspace*{-0.7cm}
\caption{Situation from the proof of Theorem~\ref{thm:threeoddcylces}.}
\label{fig:exactly-three-odd}
\end{figure}

To complete the proof we are going to show that the cycles $C$ and $D$ together with the path $Q$ in any case yield two odd cycles in addition to $C$ and $D$. For this sake we introduce some more notation. Let $z$ and $z'$ be the endvertices of the path $P$, let $P_1$ be the $x,z$-path on $C$ not passing through $z'$ and let $P_1'$ be the $x,z'$-path on $C$ not passing through $z$, see Fig.~\ref{fig:exactly-three-odd}. 
Similarly, let $P_2$ be the $y,z$-path on $D$ not passing through $z'$ and let $P_2'$ be the $y,z'$-path on $D$ not passing through $z$. Let $p_i$ be the length of $P_i$, $i\in [2]$, and let $q$ be the length of $Q$. We distinguish the following cases. 

\medskip\noindent
{\bf Case 1.} $q$ is even.\\
Assume $p_1$ and $p_2$ are both even. Then $Q\cup P_1\cup P\cup P_2'$ is an odd cycle because the length of the subpath $P\cup P_2'$ is odd. Similarly, $Q\cup P_2\cup P\cup P_1'$ is an odd cycle. 

Assume $p_1$ and $p_2$ are both odd. Then the same two cycles as in the previous subcase are odd, but this time because the length of $P-P_2'$ is even. 

In the last subcase assume that $p_1$ is odd and $p_2$ is even. (The case when $p_1$ is even and $p_2$ is odd is symmetric.) Now the extra odd cycles we are searching for are $Q\cup P_1\cup P_2$ and $Q\cup P_1'\cup P_2'$. The first cycle is clearly odd, while the second is odd because the cycle $C\cdot D$ is even and hence the length of the subpath $P_1'\cup P_2'$ is odd.      
  
\medskip\noindent
{\bf Case 2.} $q$ is odd.\\
If $p_1$ and $p_2$ are both even, or if $p_1$ and $p_2$ are both odd, then the extra odd cycles are again $Q\cup P_1\cup P_2$ and $Q\cup P_1'\cup P_2'$. 

Assume that $p_1$ is odd and $p_2$ is even. (Again the case when $p_1$ is even and $p_2$ is odd is symmetric.) Then $Q\cup P_1'\cup P\cup P_2$ are $Q\cup P_1\cup P\cup P_2$ are the required odd cycle. For instance, the latter cycle is odd because each of the subpaths $Q$, $P_1$ and $P\cup P_2$ has odd length.   
\qed

\begin{lemma}
\label{lem:intersection-path}
If $G$ is a $(3,2)$-critical graph that contains at least three odd cycles, then there exist two odd cycles whose intersection is a path with at least one edge. 
\end{lemma}
\proof
Since $G$ is $(3,2)$-critical, any two odd cycles intersect in at least two vertices. (Indeed, if odd cycles $C$ and $D$ intersect in one vertex or not at all, then there exists at least one edge $e$ which is not in $C\cup D$, and so $\es(G-e)\ge 2$, a contradiction.)
Let $C$ and $D$ be two odd cycles in $G$, and suppose that their intersection is not a path. Among vertices from $C\cap D$ between which there is no path lying in $C\cap D$, let $x$ and $y$ be chosen to be closest on $D$. Therefore, there is an $x,y$-path $P$ on $D$, which is internally disjoint with $C$, and also an $x,y$-path $Q$ on $C$, which is internally disjoint with $D$. Now, either $P\cup Q$ or $(C-Q)\cup P$ is an odd cycle; in any case the intersection of this odd cycle with $C$ is a path. 
\qed

\begin{lemma}
\label{lem:everytwointesect}
If $G$ is a $(3,2)$-critical graph that contains at least three odd cycles, then every two distinct odd cycles intersect in more than one vertex.
\end{lemma}
\proof
Suppose that $D_1$ and $D_2$ are odd cycles in $G$, which either intersect in one vertex or they are disjoint. Since the edges of $D_1\cup D_2$ induce a graph with exactly two odd cycles, and $G$ has more than two odd cycles, there must exist an edge $e$ in $G$, which is not in $D_1\cup D_2$. Now, $G-e$ still has the cycles $D_1$ and $D_2$ (having no edge in their intersection), and so $\es(G-e)\ge 2$, a contradiction.
\qed

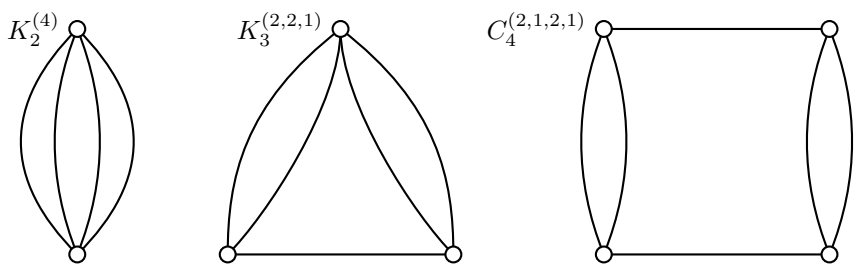
\begin{figure}[ht!]
\begin{center}
\begin{tikzpicture}[scale=1.0,style=thick]
\def\vr{3pt}
\path (0,0) coordinate (a);
\path (0,3) coordinate (b);
\path (2,0) coordinate (c);
\path (3.5,3) coordinate (d);
\path (5,0) coordinate (e);
\path (7,0) coordinate (f);
\path (7,3) coordinate (g);
\path (10,0) coordinate (h);
\path (10,3) coordinate (i);
\draw (c) -- (e); 
\draw (g) -- (i); 
\draw (f) -- (h); 
\draw (a) .. controls (-1,1) and (-1,2) .. (b);
\draw (a) .. controls (-0.4,1) and (-0.4,2) .. (b);
\draw (a) .. controls (1,1) and (1,2) .. (b);
\draw (a) .. controls (0.4,1) and (0.4,2) .. (b);
\draw (c) .. controls (2,1) and (2.2,2) .. (d);
\draw (c) .. controls (2.5,0.5) and (3.5,2) .. (d);
\draw (e) .. controls (5,1) and (4.8,2) .. (d);
\draw (e) .. controls (4.5,0.5) and (3.5,2) .. (d);
\draw (f) .. controls (6.6,1) and (6.6,2) .. (g);
\draw (f) .. controls (7.4,1) and (7.4,2) .. (g);
\draw (h) .. controls (9.6,1) and (9.6,2) .. (i);
\draw (h) .. controls (10.4,1) and (10.4,2) .. (i);
\draw (a) [fill=white] circle (\vr);
\draw (b) [fill=white] circle (\vr);
\draw (c) [fill=white] circle (\vr);
\draw (d) [fill=white] circle (\vr);
\draw (e) [fill=white] circle (\vr);
\draw (f) [fill=white] circle (\vr);
\draw (g) [fill=white] circle (\vr);
\draw (h) [fill=white] circle (\vr);
\draw (i) [fill=white] circle (\vr);
\draw [left] (b)++(-0.1,0) node {$K_2^{(4)}$};
\draw [left] (d)++(-0.1,0) node {$K_3^{(2,2,1)}$};
\draw [left] (g)++(-0.1,0) node {$C_4^{(2,1,2,1)}$};
\end{tikzpicture}
\end{center}
\caption{Multigraphs $K_2^{(4)}$, $K_3^{(2,2,1)}$, and $C_4^{(2,1,2,1)}$}
\label{fig:three-multigraphs}
\end{figure}

We denote by $\KD$ the multigraph on two vertices connected by four parallel edges, 
by $\KT$ the multigraph on three vertices two pairs of which are connected by two parallel edges and the third pair with a single edge, and by $\CF$ the multigraph on four vertices obtained from the graph $C_4$ by duplicating two of its non-consecutive edges. See Fig.~\ref{fig:three-multigraphs}.

\begin{proposition}
\label{prp:three-subdivisions}
If $G$ is a $(3,2)$-critical graph that contains at least three odd cycles, then $G$ contains as a subgraph a subdivision of one of the multigraphs $\KD$, $K_4$, $\KT$, or $\CF$. Moreover, each of the subdivisions contains at least two odd cycles. 
\end{proposition} 

\proof
By Lemma~\ref{lem:intersection-path}, there exist two odd cycles $D_1$ and $D_2$ in $G$ whose intersection is a path with at least one edge. 
We claim that there is a path $P$ connecting two distinct vertices from $D_1\cup D_2$, which is internally disjoint with $D_1\cup D_2$. Since the edges of $D_1\cup D_2$ induce only two odd cycles, that is $D_1$ and $D_2$, a third odd cycle $D_3$ must have an edge $e$, which is not in $D_1\cup D_2$. If both endvertices of $e$ are in $D_1\cup D_2$, then the claim is true, since $e$ itself induces a desired path $P$. Suppose exactly one of the endvertices of $e=xy$ is in $D_1\cup D_2$, say $x\in V(D_1\cup D_2)$. By Lemma~\ref{lem:everytwointesect}, $D_3$ intersects with $D_i$, $i\in [2]$, in at least two vertices, hence there exists a path on $D_3$ from $y$ to $D_1\cup D_2$, which confirms the claim in this case. Finally, if both $x$ and $y$ are not in $D_1\cup D_2$, by Lemma~\ref{lem:everytwointesect} we again infer that there is path on $D_3$ from $x$ to $D_1\cup D_2$ and also a path on $D_3$ from $y$ to $D_1\cup D_2$. These two paths together with the edge $e=xy$ yield a desired path $P$.

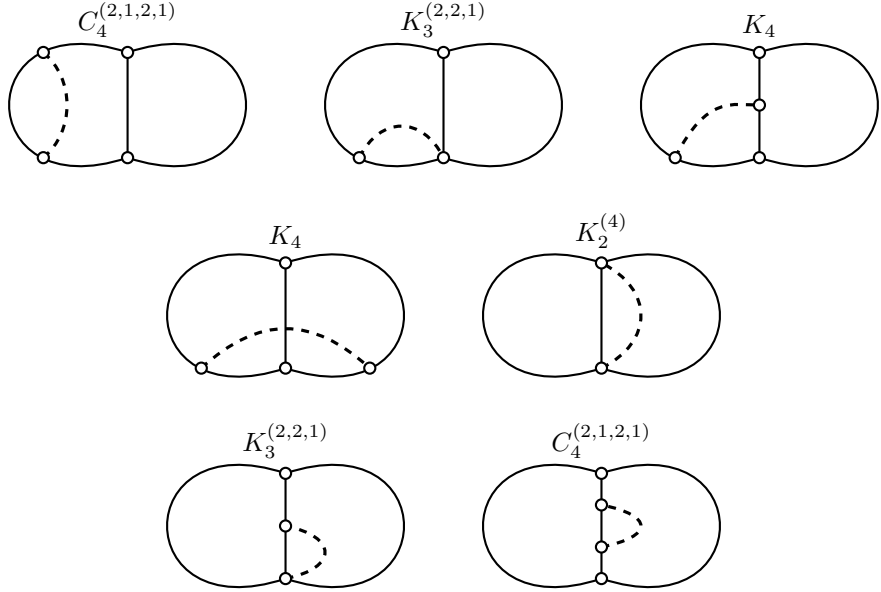
\begin{figure}[ht!]
\begin{center}
\begin{tikzpicture}[scale=0.7,style=thick]
\def\vr{3pt}
\begin{scope}
\path (0,0) coordinate (a);
\path (0,2) coordinate (b);
\path (-1.6,0) coordinate (x);
\path (-1.6,2) coordinate (y);
\draw (a) -- (b); 
\draw (a) .. controls (-3,-1) and (-3,3) .. (b);
\draw (a) .. controls (3,-1) and (3,3) .. (b);
\draw [very thick, dashed] (x) .. controls (-1,0.5) and (-1,1.5) .. (y);
\draw (a) [fill=white] circle (\vr);
\draw (b) [fill=white] circle (\vr);
\draw (x) [fill=white] circle (\vr);
\draw (y) [fill=white] circle (\vr);
\draw [above] (b)++(0,0.1) node {$C_4^{(2,1,2,1)}$};
\end{scope}
\begin{scope}[xshift = 6cm]
\path (0,0) coordinate (a);
\path (0,2) coordinate (b);
\path (-1.6,0) coordinate (x);
\draw (a) -- (b); 
\draw (a) .. controls (-3,-1) and (-3,3) .. (b);
\draw (a) .. controls (3,-1) and (3,3) .. (b);
\draw (x) [very thick, dashed] .. controls (-1.2,0.8) and (-0.4,0.8) .. (a);
\draw (a) [fill=white] circle (\vr);
\draw (b) [fill=white] circle (\vr);
\draw (x) [fill=white] circle (\vr);
\draw [above] (b)++(0,0.1) node {$K_3^{(2,2,1)}$};
\end{scope}
\begin{scope}[xshift = 12cm]
\path (0,0) coordinate (a);
\path (0,2) coordinate (b);
\path (-1.6,0) coordinate (x);
\path (0,1) coordinate (y);
\draw (a) -- (b); 
\draw (a) .. controls (-3,-1) and (-3,3) .. (b);
\draw (a) .. controls (3,-1) and (3,3) .. (b);
\draw (x) [very thick, dashed] .. controls (-1,1.2) and (-0.5,1) .. (y);
\draw (a) [fill=white] circle (\vr);
\draw (b) [fill=white] circle (\vr);
\draw (x) [fill=white] circle (\vr);
\draw (y) [fill=white] circle (\vr);
\draw [above] (b)++(0,0.1) node {$K_4$};
\end{scope}
\begin{scope}[xshift = 3cm, yshift = -4cm]
\path (0,0) coordinate (a);
\path (0,2) coordinate (b);
\path (-1.6,0) coordinate (x);
\path (1.6,0) coordinate (y);
\draw (a) -- (b); 
\draw (a) .. controls (-3,-1) and (-3,3) .. (b);
\draw (a) .. controls (3,-1) and (3,3) .. (b);
\draw (x) [very thick, dashed] .. controls (-0.5,1) and (0.5,1) .. (y);
\draw (a) [fill=white] circle (\vr);
\draw (b) [fill=white] circle (\vr);
\draw (x) [fill=white] circle (\vr);
\draw (y) [fill=white] circle (\vr);
\draw [above] (b)++(0,0.1) node {$K_4$};
\end{scope}
\begin{scope}[xshift = 9cm, yshift = -4cm]
\path (0,0) coordinate (a);
\path (0,2) coordinate (b);
\path (-1.6,0) coordinate (x);
\draw (a) -- (b); 
\draw (a) .. controls (-3,-1) and (-3,3) .. (b);
\draw (a) .. controls (3,-1) and (3,3) .. (b);
\draw [very thick, dashed] (a) .. controls (1,0.5) and (1,1.5) .. (b);

\draw (a) [fill=white] circle (\vr);
\draw (b) [fill=white] circle (\vr);
\draw [above] (b)++(0,0.1) node {$K_2^{(4)}$};
\end{scope}
\begin{scope}[xshift = 3cm, yshift = -8cm]
\path (0,0) coordinate (a);
\path (0,2) coordinate (b);
\path (0,1) coordinate (x);
\draw (a) -- (b); 
\draw (a) .. controls (-3,-1) and (-3,3) .. (b);
\draw (a) .. controls (3,-1) and (3,3) .. (b);
\draw (a) [very thick, dashed] .. controls (1,0.2) and (1,0.8) .. (x);
\draw (a) [fill=white] circle (\vr);
\draw (b) [fill=white] circle (\vr);
\draw (x) [fill=white] circle (\vr);
\draw [above] (b)++(0,0.1) node {$K_3^{(2,2,1)}$};
\end{scope}
\begin{scope}[xshift = 9cm, yshift = -8cm]
\path (0,0) coordinate (a);
\path (0,2) coordinate (b);
\path (0,0.6) coordinate (x);
\path (0,1.4) coordinate (y);
\draw (a) -- (b); 
\draw (a) .. controls (-3,-1) and (-3,3) .. (b);
\draw (a) .. controls (3,-1) and (3,3) .. (b);
\draw (x) [very thick, dashed] .. controls (1,0.8) and (1,1.2) .. (y);
\draw (a) [fill=white] circle (\vr);
\draw (b) [fill=white] circle (\vr);
\draw (x) [fill=white] circle (\vr);
\draw (y) [fill=white] circle (\vr);
\draw [above] (b)++(0,0.1) node {$C_4^{(2,1,2,1)}$};
\end{scope}
\end{tikzpicture}
\end{center}
\vspace*{-0.7cm}
\caption{Situations leading to subdivisions of $K_2^{(4)}$, of $K_3^{(2,2,1)}$, of $C_4^{(2,1,2,1)}$, and of $K_4$}
\label{fig:possible-subdivisions}
\end{figure}

All the described possibilities for the position of $P$ in $D_1\cup D_2$ are schematically shown in Fig.~\ref{fig:possible-subdivisions}.
We infer that $G$ contains as a subgraph a subdivision of one of the following graphs: $\KD$, $K_4$, $\KT$, or $\CF$. Since the odd cycles $D_1$ and $D_2$ are a part of this subdivision, the last sentence of the statement of the proposition is also clear. 
\qed

\begin{theorem}
\label{thm:fouroddcycles}
The $(3,2)$-critical graphs, which contain exactly four odd cycles are precisely the graphs of the families $\cal C$ and $\cal D$.
\end{theorem}
\proof
Let $G$ be a $(3,2)$-critical graph, which contains exactly four odd cycles $D_1$, $D_2$ , $D_3$, and $D_4$. We start with the following claim.
\medskip

\noindent {\bf Claim.} Every edge of $G$ is contained in at least two odd cycles.

\smallskip

\begin{figure}[ht!]
\begin{center}
\begin{tikzpicture}[scale=0.9,style=thick]
\def\vr{3pt}


\path (0,0) coordinate (a);
\path (2,0) coordinate (b);
\path (3,0) coordinate (x);
\path (-1,0) coordinate (y);
\draw (a) -- (b); \draw (a) -- (y); \draw (b) -- (x); 
\draw (y) .. controls (-3,0) and (0,-5) .. (b);
\draw (x) .. controls (5,0) and (2,-5) .. (a);
\draw (x) .. controls (2.5,2) and (-0.5,2) .. (y);
\draw (a) [fill=white] circle (\vr);
\draw (b) [fill=white] circle (\vr);
\draw (1,-0.2) node {$P$};
\draw (-0.4,-1.2) node {$D_1$};
\draw (2.4,-1.2) node {$D_2$};
\draw (1,0.7) node {$D_3$};

\end{tikzpicture}
\end{center}
\vspace*{-3cm}
\caption{Case B from the proof of Claim in the proof of Theorem~\ref{thm:fouroddcycles}}
\label{fig:two-cases}
\end{figure}
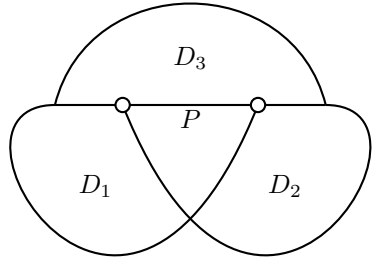

\noindent {\bf Proof} (of Claim).  Suppose $e$ is an edge which is contained only in one of the cycles, say $D_4$. Hence, since $\es(G-e)=1$, there exists an edge $f$ in $D_1\cap D_2\cap D_3$. By Lemma~\ref{lem:general}, the intersection of every two odd cycles in $G-e$ is a path. Then, $D_1\cap D_2\cap D_3$ is also a path (containing $f$), and let us denote it by $P$. Note that for two of the cycles among $D_1,D_2,D_3$ their intersection is the path $P$. Without loss of generality, let $D_1\cap D_2=P$. There are two cases: either $D_1$ (resp. $D_2$) has $D_1\cap D_3=P$ (resp. $D_2\cap D_3=P$), or, $(D_1\cap D_3)-P\ne\emptyset$ and $(D_2\cap D_3)-P\neq \emptyset$; see Fig.~\ref{fig:two-cases} for the second case.

\medskip

\noindent {\bf Case A.} $D_1\cap D_3=P$.\\
Note that $E(D_4\cap P)=\emptyset$, that is, $E(D_4)\cap E(D_1\cap D_2)=\emptyset=E(D_4)\cap E(D_1\cap D_3)$. Since $D_4$ contains an edge which does not lie in $D_1\cup D_2\cup D_3$, there is an edge $g$ in the (even) cycle $(D_2\cup D_3)-(D_2\cap D_3)$, which is not in $D_4$. Note that $g$ also does not lie in $D_1$.  Now, we claim that $\es(G-g)\ge 2$. Indeed, if an edge in $P$ is removed from $G-g$, then $D_4$ is still an odd cycle in the resulting graph. On the other hand, if an edge in $D_1-P$ is removed from $G-g$, then either $D_2$ or $D_3$ remains in the resulting graph (depending on whether $g$ is in $D_2$ or $D_3$). Otherwise, $D_1$ remains in the resulting graph obtained by the removal of two edges, which is the final contradiction in this case.

\medskip

\noindent {\bf Case B.} $(D_1\cap D_3)-P\ne\emptyset$ and $(D_2\cap D_3)-P\neq \emptyset$.\\
Note that $((D_1\cup D_2\cup D_3)-((D_1\cap D_3)\cup(D_2\cap D_3)))\cup P$ induces an odd cycle. Since $G$ has exactly four odd cycles, this cycles must be $D_4$. This is in a contradiction to the assumption that $e$ is an edge of $D_4$, not contained in the other three cycles. 

\medskip

Both cases lead to a contradiction, hence every edge of $G$ indeed lies in at least two odd cycles, as claimed. \qqed

\medskip

By Proposition~\ref{prp:three-subdivisions}, since $G$ contains at least three odd cycles, it contains a subgraph, which is a subdivision of one of the multigraphs $\KD$, $K_4$, $\KT$, or $\CF$. In all the cases let $G'$ denote the respective subgraph, a subdivision of the corresponding multigraph. 

\medskip\noindent
{\bf Case 1.} $G$ contains a subdivision of $\KD$.\\
Let $x$ and $y$ be the vertices in $G$, whose degree in $G'$ is $4$, and let $P_i$, $i\in [4]$, be the corresponding $x,y$-paths in $G'$, having lengths $p_i$, respectively. If exactly two of the lengths $p_i$ are odd, then $G'$ is clearly a $(3,2)$-critical graph. If $E(G)-E(G')\ne\emptyset$, then the removal of an edge in $E(G)-E(G')$ yields a graph with $\es\ge 2$. Therefore, $G=G'$, and $G$ is in $\cal C$. If an integer $p_j$ is of different parity than the other three, then $G'$ contains exactly three odd cycles, all of which contain the path $P_j$. Hence $\es(G')=1$. By Claim, every edge of $G$ lies in at least two odd cycles, thus the fourth odd cycle of $G$ contains all the edges of $E(G')-E(P_j)$, which is not possible. Therefore in this case, $G$ is not $(3,2)$-critical. In the other two cases for parities of integers $p_i$, $G'$ is bipartite, which is not a relevant case according to Proposition~\ref{prp:three-subdivisions}. 

\medskip\noindent
{\bf Case 2.} $G$ contains a subdivision of $K_4$.\\
In this case, $G'$ has six paths $P_i$ with lengths $p_i$ that connect vertices of degree $3$ in $G'$, and there are seven cycles, four of which come from the triangles of $K_4$ and three of which come from the squares of $K_4$. If $G'$ is in $\cal D$, then note that $G'$ has exactly four odd cycles, and every edge of $G'$ lies in exactly two odd cycles. As in Case 1, we infer that $G=G'$, and so $G$ is in $\cal D$. In other cases of parities of integers $p_i$, $G'$ contains exactly four odd cycles, all of which pass through the (common) path $P_i$ of odd length. This implies $\es(G')=1$, and since $G$ has no more odd cycles, $\es(G)=1$. (The case when all parities are even is not relevant according to Proposition~\ref{prp:three-subdivisions}.)

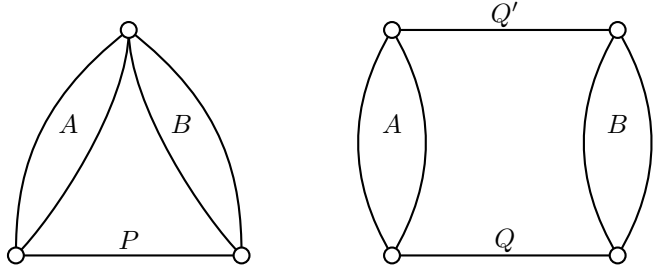
\begin{figure}[ht!]
\begin{center}
\begin{tikzpicture}[scale=1.0,style=thick]
\def\vr{3pt}
\path (2,0) coordinate (c);
\path (3.5,3) coordinate (d);
\path (5,0) coordinate (e);
\path (7,0) coordinate (f);
\path (7,3) coordinate (g);
\path (10,0) coordinate (h);
\path (10,3) coordinate (i);
\draw (c) -- (e); 
\draw (g) -- (i); 
\draw (f) -- (h); 
\draw (c) .. controls (2,1) and (2.2,2) .. (d);
\draw (c) .. controls (2.5,0.5) and (3.5,2) .. (d);
\draw (e) .. controls (5,1) and (4.8,2) .. (d);
\draw (e) .. controls (4.5,0.5) and (3.5,2) .. (d);
\draw (f) .. controls (6.4,1) and (6.4,2) .. (g);
\draw (f) .. controls (7.6,1) and (7.6,2) .. (g);
\draw (h) .. controls (9.4,1) and (9.4,2) .. (i);
\draw (h) .. controls (10.6,1) and (10.6,2) .. (i);
\draw (c) [fill=white] circle (\vr);
\draw (d) [fill=white] circle (\vr);
\draw (e) [fill=white] circle (\vr);
\draw (f) [fill=white] circle (\vr);
\draw (g) [fill=white] circle (\vr);
\draw (h) [fill=white] circle (\vr);
\draw (i) [fill=white] circle (\vr);
\draw (3.5,0.2) node {$P$};
\draw (2.7,1.75) node {$A$};
\draw (4.2,1.75) node {$B$};
\draw (8.5,0.2) node {$Q$};
\draw (8.5,3.2) node {$Q'$};
\draw (7,1.75) node {$A$};
\draw (10,1.75) node {$B$};
\end{tikzpicture}
\end{center}
\caption{For the case analysis}
\label{fig:last-case-analysis}
\end{figure}

\medskip\noindent
{\bf Case 3.} $G$ contains a subdivision of $\KT$.\\
In this case, $G'$ consists of five (subdivision) paths, two pairs of which form cycles $A$ and $B$, and $P$ is the path not in $A\cup B$; see the left graph in Fig.~\ref{fig:last-case-analysis}. 

If both cycles $A$ and $B$ are odd, then by Lemma~\ref{lem:everytwointesect}, $G'$ is not $(3,2)$-critical. In fact, $\es(G')\ge 3$, and so $\es(G)\ge 3$. If both $A$ and $B$ are even cycles (and $G'$ is not bipartite), then $G'$ has four odd cycles all of which have path $P$ in common. Hence, $\es(G')=1$, and since $G$ has no other odd cycles, also $\es(G)=1$. If one of the two cycles, say $A$, is odd, and $B$ is even, then $G'$ contains exactly three odd cycles, which all have a path $R$ from $A$ in common. Hence $\es(G')=1$. Since every edge of $G$ lies in at least two odd cycles (by Claim), the fourth odd cycle of $G$ contains all the edges of $E(G')-E(R)$, which is not possible. Therefore also in this case, $G$ is not $(3,2)$-critical.

\medskip\noindent
{\bf Case 4.} $G$ contains a subdivision of $\CF$.\\
In this case, $G'$ consists of six paths, two pairs of which form cycles $A$ and $B$, and  call the other two paths by $Q$ and $Q'$; see the right graph in Fig.~\ref{fig:last-case-analysis}. If both cycles $A$ and $B$ are odd, we infer as in the previous case that $G$ is not $(3,2)$-critical, using Lemma~\ref{lem:everytwointesect}. If $A$ and $B$ are both even cycles, then $G'$ has four odd cycles, which all go through $Q$ or $Q'$. Without loss of generality, assume that all four cycles pass through $Q$. Hence $\es(G')=1$. Since every edge of $G$ lies in at least two odd cycles, the fourth odd cycle of $G$ contains all the edges of $E(G')-E(Q)$, which is not possible. For the final case, assume that $A$ is odd and $B$ is even (the reversed case is symmetric). Then $G'$ contains exactly three odd cycles, which are passing through one of the paths in $A$. Hence $\es(G')=1$, and we again infer that the fourth odd cycle of $G$ contains all edges of $G'$ except those of one of the paths in $A$, which is not possible. 
\qed

Theorem~\ref{thm:main} now follows by combining Theorems~\ref{thm:two-odd-cycles}, \ref{thm:threeoddcylces}, and~\ref{thm:fouroddcycles}.

\section{On the complexity of recognizing $(k,2)$-critical graphs}
\label{sec:complexity}

In this section, we investigate the computational complexity of the recognition of $(k,2)$-critical graphs for $k\ge 4$. It is easy to see that $(3,2)$-critical graphs can be efficiently recognized. Indeed, for every edge $e$ of a given graph $G$ one has to verify whether $G-e$ is not bipartite and whether there exists an edge $f$ such that $G-\{e,f\}$ is bipartite; clearly, a BFS search can be used to check (non)bipartiteness, which yields a polynomial algorithm for recognizing $(3,2)$-critical graphs. On the other hand, the subsequent result in this section indicates that the recognition of $(k,2)$-critical graphs when $k\ge 4$ is likely to be computationally hard. 

Given an integer $k$, $k\ge 3$, let $\KKD$ denote the class of graphs that are $(k,2)$-critical, and let $\KCH$ denote the class of graphs $G$ with $\chi(G)=k$ that are critical for the chromatic number. That is, $G$ is in $\KCH$ when $\chi(G)=k$ and $\chi(H)<k$ for each proper subgraph $H$ of $G$. Since we consider only graphs $G$ with no isolated vertices, the latter condition is equivalent to the statement that $\chi(G-e)=k-1$ for every edge $e$ in $G$. For more information on graphs critical for the chromatic number see Chapter 5 in the book~\cite{jensen-1995}, recent papers~\cite{kostochka-2018, postle-2018}, and references therein. 

In the next result we show that the problem of deciding whether a graph is in $\KCH$ can be reduced in polynomial time to the problem of deciding whether a graph is in $\KKD$. In particular, the existence of a polynomial algorithm for recognizing  the graphs from $\KKD$ implies that there exists a polynomial algorithm for recognizing the graphs is in $\KCH$. Since the latter seems very unlikely, we believe that recognizing the graphs in $\KKD$ is computationally hard. 

Let $G$ be a connected graph and $u\in V(G)$. We denote with 
$G_u \bowtie K_k$ the graph obtained from $G$ and $K = K_k$ by identifying $u$ with a vertex of $K$. Let further $K_u$ be the $k$-clique of $G_u \bowtie K_k$ obtained in the construction from $K$.  

\begin{theorem}
Let $G$ be a graph and let $u\in V(G)$. If $k\ge 3$, then $G\in\KCH$ if and only if $G_u \bowtie K_k\in \KKD$. 
\end{theorem}

\proof
Let $G\in \KCH$.  Clearly, $\chi(G_u \bowtie K_k)=k$. To prove that $G_u \bowtie K_k\in \KKD$, consider two cases for an edge $e$ in $G_u \bowtie K_k$.  If $e$ is in $K_u$, then clearly $\chi((G_u \bowtie K_k)-e)=k$, since $G$ is a subgraph of $(G_u \bowtie K_k)-e$. In addition, for an arbitrary edge $f\in E(G)$, we have $\chi((G_u \bowtie K_k)-\{e,f\})=k-1$, and so $\es(G-e)=1$. The other case for an edge $e$ (that is, $e\in E(G)$) is similar. Indeed, we have $\chi((G_u \bowtie K_k)-e)=k$, yet for an arbitrary edge $f$ from $K_u$ we have $\chi((G_u \bowtie K_k)-\{e,f\})=k-1$. Thus, $G_u \bowtie K_k\in \KKD$.

For the reversed implication, let $G_u \bowtie K_k\in \KKD$. Since $\chi((G_u \bowtie K_k)-e)=k$ for an edge $e\in E(K_u)$, this implies that $\chi(G)=k$. To prove that $G$ is in $\KCH$, let $e$ be an arbitrary edge in $E(G)$. Note that $\chi((G_u \bowtie K_k)-e)=k$. Since $\es((G_u \bowtie K_k)-e)=1$, there exists an edge $f\in E(G_u \bowtie K_k)$ such that $\chi((G_u \bowtie K_k)-\{e,f\})=k-1$. This edge $f$ cannot lie in $E(G)$ for otherwise $(G_u \bowtie K_k)-\{e,f\}$ contains the $k$-clique $K_u$ as a subgraph, and so the chromatic number of this graph is $k$. So $f\in E(K_u)$ and hence $\chi(G-e)=k-1$.
\qed

\section{Concluding remarks on $(3,2)$-critical graphs}

We could find only one more family of $(3,2)$-critical graphs. They can be described as specific subdivisions of a multigraph obtained from a cycle whose edges are duplicated. More precisely, a graph from $\cal E$ is obtained from the disjoint union of $k$ even cycles $C_{2n_1},\ldots, C_{2n_k}$ as follows. For each $i\in [k]$, let $x_i$ and $y_i$ be any two distinct vertices of $C_{2n_i}$, where we only require that $$\sum_{i=1}^k d_{C_{2n_i}}(x_i,y_i)$$ is odd. A graph $G\in \cal E$ is obtained by identifying $y_i$ with $x_{i+1}$ for $i\in [k-1]$, and identifying $y_k$ with $x_1$. It is easy to see that graphs $G\in \cal E$ are $(3,2)$-critical. We pose the following problem for which we suspect it has a positive answer.

\begin{problem}
\label{prob}
Is it true that $\cal A\cup \cal B\cup \cal C\cup \cal D\cup \cal E$ is the family of $(3,2)$-critical graphs (without isolated vertices)? 
\end{problem}

Solving Problem~\ref{prob} remains a challenge. One way how to approach it is by using Proposition~\ref{prp:three-subdivisions}, which shows that a $(3,2)$-critical graph with at least three odd cycles contains one of the four subdivisions described in the proposition as a subgraph. Now, assuming that a $(3,2)$-critical graph $G$ has at least five odd cycles, one should examine each of the four possible subdivisions appearing as a proper subgraph $H$ of $G$, and each possibility of which are the odd cycles in these subdivisions. Cases of $H$ being a subdivision of $\KD$ and $K_4$ should yield a contradiction, while $\KT$ and $\CF$ should either yield a contradiction or that $G$ belongs to $\cal E$. The described approach is probably very technical, hence applying some related work from the literature in an efficent way would be desirable. Among the known results about smallest sets of edges that hit all odd cycles we encountered papers of Berge and Reed~\cite{ber-2000} and Kr\'{a}l and Voss~\cite{kral-2004} whose main concern are planar graphs, and cannot be applied for this purpose.

\section*{Acknowledgments}

The third author thankfully acknowledges initial discussions on the problem with Ross Kang.  We also thank Jason Brown and Tommy Jensen for discussions on the recognition complexity of chromatic critical graphs. The financial support from the Slovenian Research Agency (research core funding P1-0297, projects J1-9109, J1-1693, and N1-0095) is acknowledged. The third author acknowledges the financial support from Yazd University research affairs as Post-doc research project.

\end{document}